\documentclass{article}
\usepackage[latin1]{inputenc}
\usepackage{latexsym}
\usepackage{a4wide}
\usepackage{amscd}
\usepackage{graphics}
\usepackage{amsmath}
\usepackage{amssymb}
\usepackage{mathrsfs}
\input xy
\xyoption{all}

\newcommand{\Z}{\mathbb{Z}}                     
\newcommand{\R}{\mathbb{R}}                     
\newcommand{\C}{\mathbb{C}}                     
\newcommand{\proof}{{\sl Proof.}\hspace{5pt}}   
\newcommand{\finedim}{\hfill $\Box$}            
\newcommand{\Gr}{\mathrm{Gr}}
\newcommand{\sgn}{\mathrm{sign}}

\newcommand{\Spl}{\mathrm{Sp}}
\newcommand{\spl}{\mathfrak{sp}}

\newcommand{\diag}{\mathrm{diag}}

\numberwithin{equation}{section}

\newtheorem{thm}{\sc Theorem}[section]      
\newtheorem{cor}[thm]{\sc Corollary}        
\newtheorem{lem}[thm]{\sc Lemma}            
\newtheorem{prop}[thm]{\sc Proposition}     
\newtheorem{defin}[thm]{\sc Definition}     
\newtheorem{remark}[thm]{\sc Remark}        

\title{Maslov index of Hamiltonian systems}

\author{Alessandro Portaluri}

\date{\today}

\begin{document}
\maketitle


\begin{abstract}
The aim of this paper is to give an explicit formula in order to
compute the Maslov index of the fundamental solution of a linear
autonomous Hamiltonian systems in terms of the Conley-Zehnder index
and the time one flow.
\end{abstract}
\section*{Introduction}\label{sec:intro}

The Maslov index is a semi-integer homotopy invariant of paths $l$
of Lagrangian subspaces of a symplectic vector space $(V, \omega)$
which gives the algebraic counts of non transverse intersections of
the family $\{l(t)\}_{t \in [0,1]}$ with a given Lagrangian subspace
$l_{*}$. To be more precise, let us denote by
$\Lambda(V):=\Lambda(V, \omega)$ the set of all Lagrangian subspaces
of the symplectic space $V$ and let $\Sigma(l_*)  = \{ l \in
\Lambda: l \cap l_* \not= (0)\}$ be the {\em train\/} or the {\em
Maslov cycle\/} of $l_*$. Then, it can be proven that $\Sigma(l_*)$
is a co-oriented one codimensional algebraic subvariety of the
Lagrangian Grassmannian $\Lambda(V)$ and the Maslov index counts
algebraically the number of intersections of $l$ with $\Sigma(l_*)$.
This is the basic invariant out of which many others are defined.
For example, if $\phi \colon [a,b] \to \Spl(V)$ is a path of
symplectic automorphisms of $V$ and $l_*$ is a fixed Lagrangian
subspace, then the Maslov index of $\phi$ is by definition the
number of intersections of the path $[a,b]\ni t \mapsto
\phi_{t}(l_*)\in \Lambda(V)$ with the train of $l_*$. The aim of
this paper is to explicitly compute the Maslov index of the
fundamental solution associated to
\[
w'(x) =  H  w(x)
\]
where $H$ is a (real) constant Hamiltonian matrix. The idea in order
to perform our computation is to relate the Maslov index with the
Conley-Zehnder index and then to compute a correction term arising
in this way which is expressed in terms of an invariant of a triple
of Lagrangian subspaces also known as Kashiwara index.

We stress the fact that the result is not new and it was already
proven in \cite{MarPicTau01}. However the contribution of this paper
is to provide a different and we hope a simpler proof of this
formula.

\section{Linear preliminaries}\label{sec:symplecticgroup}

The purpose of this section is to recollect some well-known facts
about the geometry of the Lagrangian Grassmannian and the Maslov
index needed in our computation.  For further details see for
instance \cite{Abbo1, Dui76, GiaPicPor, MarPicTau01,RobSal}.
\begin{defin}
Let $V$ be a finite dimensional real vector space. A symplectic form
$\omega$ is a non degenerate anti-symmetric bilinear form on $V$. A
symplectic vector space is a pair $(V, \omega)$.
\end{defin}
The archetypical example of symplectic space is $(\R^{2n},
\omega_0)$ where the symplectic structure $\omega_0$ is defined as
follows. Given the splitting  $\R^{2n}= \R^n \oplus \R^n$ and the
scalar product of $\R^n$  $\langle \cdot, \cdot \rangle$ then for
each $z_k=(x_k, y_k) \in \R^n \oplus \R^n$ for $k=1,2$ we have
\[
\omega_0(z_1, z_2) =  \langle x_1, y_2 \rangle- \langle x_2, y_1
\rangle.
\]
This symplectic structure $\omega_0$ can be represented against the
scalar product by setting $ \omega_0 (z_1, z_2) = \langle J z_1, z_2
\rangle$ for all $z_i \in \R^{2n}$ with $i=1,2$ where we denoted by
$J$ the {\em standard complex structure\/} of $\R^{2n}$ which can be
written with respect to the canonical basis of $\R^{2n}$ as
\begin{equation}\label{eq:strutcomplessa}
J \, = \, \begin{pmatrix}
0 & -I_n \\
I_n & 0
\end{pmatrix}
\end{equation}
where $I_n$ is the $n$ by $n$ identity matrix. Given a linear
subspace of the symplectic vector space $(V, \omega)$, we define the
orthogonal of $W$ with respect to the symplectic form $\omega $ as
the linear subspace $W^{\sharp}$ given by $ W^{\sharp}\,=\,\{v \in V
|\quad \omega(u,v)=0\   \forall u\, \in W\}.$
\begin{defin} Let $W$ be a linear subspace of $V$. Then
\begin{enumerate}
\item[(i)] $W$ is {\em isotropic\/} if $W\subset W^{\sharp}$;
\item[(ii)] $W$ is {\em symplectic\/} if $ W^{\sharp}\cap W= {0}$.
\item[(iii)] $W$ is {\em Lagrangian\/} if $W = W^\sharp.$
\end{enumerate}
\end{defin}
\begin{defin}
Let $(V_1, \omega_1)$ and $(V_2, \omega_2)$ be symplectic vector
spaces. A symplectic isomorphism from $(V_1, \omega_1)$ to $(V_2,
\omega_2)$ is a bijective linear map $\varphi \colon V_1 \to V_2$
such that $\varphi^* \omega_2 = \omega_1$, meaning that
\[
\omega_2(\varphi(u), \varphi(v))= \omega_1(u,v), \ \ \forall \, u, v
\in V_1.
\]
In the case $(V_1, \omega_1)=(V_2, \omega_2)$, $\varphi$ is called a
symplectic automorphism or symplectomorphism.
\end{defin}
The matrices which correspond to symplectic automorphisms of the
standard symplectic space $(\R^{2n}, \omega_0)$ are called
symplectic and they are characterized by the equation
\[
A^T JA =J
\]
where $A^T$ denotes the adjoint of $A$. The set of all symplectic
automorphisms of $(V, \omega)$ forms a group, denoted by $\Spl(V,
\omega)$. The set of the symplectic matrices is a Lie group, denoted
by $\Spl(2n)$. Since each symplectic vector space of dimension $2n$
is symplectically isomorphic to $(\R^{2n}, \omega_0)$, then $\Spl(V,
\omega)$ is isomorphic to $\Spl(2n)$. The Lie algebra of $\Spl(2n)$
is
\[\spl (2n):= \{H \in L(2n)| H^TJ + JH=0\}
\]
where $L(2n)$ is the vector space of all real matrices of order
$2n$. The matrices in $\spl(2n)$ are called {\em infinitesimally
symplectic\/} or {\em Hamiltonian\/}.

\subsection{The Krein signature on $\Spl(2n)$}
Following the argument given in \cite[Chapter 1, Section 1.3]{Abbo1}
we briefly recall the definition of Krein signature of the
eigenvalues of a symplectic matrix.

In order to define the Krein signature of a symplectic matrix $A$ we
shall consider $A$ as acting on $\C^{2n}$ in the usual way
\[
A(\xi + i \eta) := A \xi + i A \eta, \qquad \forall \, \xi, \eta \in
\R^{2n}
\]
and we define the Hermitian form $g(\xi,\eta) := \langle G \xi, \eta
\rangle$ where  $G:=-i J$. The complex symplectic group $\Spl(2n,
\C)$ consists of the complex matrices $A$ such that
\begin{equation}\label{eq:sullacomplessadisp}
A^* GA = G
\end{equation}
where as usually $ A^*=\bar A^T$ denotes the transposed conjugate of
$A$.

\begin{defin}\label{def:signaturadokrein}
Let $\lambda$ be an eigenvalue on the unit circle of a complex
symplectic matrix. The {\em Krein signature\/} of $\lambda$ is the
signature of the restriction of the Hermitian form $g$ to the
generalized eigenspace $E_\lambda$.
\end{defin}
If the real symplectic matrix $A$ has an eigenvalue $\lambda$ on the
unit circle of Krein signature $(p,q)$, it is often convenient to
say that $A$ has $p+q$ eigenvalues $\lambda$, and that $p$ of them
are Krein-positive and $q$ which are Krein-negative.

Let $A$ be a semisimple symplectic matrix, meaning that the
algebraic and geometric multiplicity of its eigenvalues coincides.
\begin{defin}\label{def: normalisemisemplicisimplectiic3}
We say that $A$ is in {\em normal form\/} if $A = A_1 \oplus \dots
\oplus A_p$, where $A_i$ has one of the forms listed below:
\begin{itemize}
\item[(i)]
\begin{equation*} A_1 =\begin{pmatrix}
\cos \alpha  & -\sin \alpha   \\
\sin \alpha  & \cos \alpha
\end{pmatrix},\qquad \textrm{for} \ \ \  \alpha \in \R.
\end{equation*}
\item[(ii)]
\begin{equation*}
A_2=\begin{pmatrix} \mu & 0  \\ 0 & \mu^{-1}
\end{pmatrix}\qquad \textrm{for} \quad \mu \in \R, \  \textrm{and}\ \
|\mu|>1.
\end{equation*}
\item[(iii)]
\begin{equation*}
A_3 =\begin{pmatrix} \lambda \cos \alpha &-\lambda \sin \alpha&
0& 0  \\
\lambda \sin \alpha & \lambda\cos \alpha & 0 & 0  \\
0 & 0 &\lambda^{-1} \cos \alpha & -\lambda^{-1} \sin \alpha \\
0 & 0& \lambda^{-1} \sin \alpha & \lambda^{- 1}\cos \alpha
\end{pmatrix},
\end{equation*}
for $\alpha \in \R \backslash\pi \Z, \mu \in \R, | \mu|>1$.
\end{itemize}
\end{defin}

\subsection{The Maslov index} In this section we define the Maslov
index for Lagrangian and symplectic paths. Our basic reference is
\cite{RobSal}. Given a symplectic space $(V, \omega)$ we denote by
$\Lambda(V, \omega)$ the set of all Lagrangian subspaces. Let
$L_0\in \Lambda(V, \omega)$ be fixed and for all $k=0,1,\ldots,n$ we
set
\[\Lambda_k(L_0)=\big\{L\in\Lambda:\mathrm{dim}(L\cap
L_0)=k\big\}\quad\text{and}\quad\Sigma(L_0)=\cup_{k=1}^n\Lambda_k(L_0).\]
It can be proven that each stratum $\Lambda_k(L_0)$ is connected of
codimension $\frac12k(k+1)$ in $\Lambda$.

Let $l \colon [a,b] \to \Lambda$ be a $C^1$-curve of Lagrangian
subspaces. We say that $l$ has a {\em crossing\/} with the {\em
train\/} $\Sigma(L_0)$ of $L_0$ at the instant $t=t_0$ if
$l(t_0)\in\Sigma(L_0)$. At each non transverse crossing time $t_0
\in [a,b]$ we define the {\em crossing form\/} $\Gamma$ as the
quadratic form
\[
\Gamma (l, L_0, t_0)\, =\, l'(t_0)|_{l(t_0)\cap L_0}
\]
and we say that a crossing $t$ is called {\em regular\/} if the
crossing form is nonsingular. It is called {\em simple\/} if it is
regular and in addition $l(t_0) \in \Lambda_1(L_0)$.
\begin{defin}\label{def:RobSalindex}
Let $l \colon [a,b] \to \Lambda$ be a smooth curve having only {\em
regular crossings\/} we define the {\em Maslov index\/}
\begin{equation}\label{eq: MaslovdiRobSal}
\mu(l, L_0)\, := \,\frac{1}{2}\sgn\, \Gamma (l, L_0, a)\, +
\,\sum_{t\in\left]a,b\right[}\sgn\, \Gamma (l, L_0, t)\, +
\,\frac{1}{2}\sgn\, \Gamma (l, L_0, b)
\end{equation}
where the summation runs over all crossings $t$.
\end{defin}
For the properties of this number we refer to \cite{RobSal}.
Now let $\psi \colon [a,b] \to \Spl(2n)$ be a continuous path of
symplectic matrices and $L \in \Lambda(n)$ where $\Lambda(n)$
denotes the set of all Lagrangian subspaces of the symplectic space
$(\R^{2n},  \omega_0)$. Then we define the {\em Maslov index\/} of
the $\psi$ as
\[
\mu_{L}(\psi)\, := \, \mu(\psi L, L).
\]
Given the vertical Lagrangian subspace $L_0 = \{0\}\oplus \R^n$ and
assuming that $\psi$ has the block decomposition
\begin{equation}\label{eq: blockdecomposition}
\psi(t)=\begin{pmatrix} a(t) & b(t) \\ b(t) & d(t)
\end{pmatrix},
\end{equation}
then the crossing form of the path of Lagrangian subspaces $\psi
L_0$ at the crossing instant $t=t_0$ is the quadratic form
$\Gamma(\psi, t_0)\colon \ker\, b(t_0) \to \R$ given by
\begin{equation}\label{eq: calcolocrossingsimplettici}
\Gamma(\psi, t_0)(v) \, = \, - \langle d(t_0)v , b'(t_0) v \rangle
\end{equation}
where $b(t_0)$ and $d(t_0)$ are the block matrices defined in
\eqref{eq: blockdecomposition}. Lemma below will be crucial in our
final computation.
\begin{lem}\label{thm:lagrangianiprodotto}
Consider the symplectic vector space $\R^{2n} \times \R^{2n}$
equipped by the symplectic form $\overline\omega = - \omega_0 \times
\omega_0$. Then
\begin{equation}\label{eq:unadellaultime}
\mu(\psi L, L_1) \,=\, \mu( \Gr(\psi),L \times L_1)
\end{equation}
where $\Gr$ denotes the graph and where $L, L_1 \in \Lambda(n)$.
\end{lem}
\proof For the proof of this result see \cite[Theorem
3.2]{RobSal}.\finedim
\begin{defin}\label{def:conleyzehnder}
Given a continuous path of symplectic matrices $\psi$, we define the
{\em Conley-Zehnder index\/} $\mu_{CZ}(\psi)$ as
\[
\mu_{CZ}(\psi):= \mu( \Gr(\psi),\Delta)
\]
where $\Delta \subset \R^{2n}\times \R^{2n}$ denotes the diagonal in
the product space.
\end{defin}
Let us consider the path
\begin{equation}\label{eq:dlpirmotipoformnorm}
\psi_1(x) =\begin{pmatrix} \cos \alpha x & -\sin \alpha x  \\ \sin
\alpha x & \cos \alpha x
\end{pmatrix}.
\end{equation}
Given the Lagrangian  $L_0= \{0\} \oplus \R$ in $\R^2$ let us
consider the path of Lagrangian subspaces of $\R^2$ given by $l_1
:=\psi_1 \,L_0$. It is easy to check that the crossing points are of
the form $x \in \pi\Z/\alpha$ so by formula \eqref{eq:
calcolocrossingsimplettici} if $x_0$ is a crossing instant then we
have:
\[
\Gamma(\psi_1, x_0)(k)= \alpha k\  \cos(\alpha x_0)k  \cos(\alpha
x_0)= \alpha k^2 \cos^2(\alpha x_0).
\]
Thus if $\alpha \not=0$ then we have
\begin{equation*}
\sgn \,\Gamma(\psi_1, x_0)\, = \, \left\{
\begin{array}{ll} 1 & \text{if} \ \alpha >0  \\
-1 & \text{if}\ \alpha < 0.
\end{array}\right.
\end{equation*}
Summing up we have
\begin{lem}\label{thm: Maslov1normale}
Let $\psi_1:[0,1] \to \Spl (2)$ be the path of symplectic matrices
given in \eqref{eq:dlpirmotipoformnorm}. Then the Maslov index is
given by
\begin{enumerate}
\item {\em non transverse end-point}
\begin{equation*}
 \mu(\psi_1)\,=\,\frac{\alpha}{\pi}
\end{equation*}
\item {\em transverse end-point}
\begin{equation*}
 \mu(\psi_1)\,=\,
\left[\frac{\alpha}{\pi}\right]
 \, + \, \frac{1}{2}
\end{equation*}
\end{enumerate}
where we have denoted by $[\cdot]$ the integer part.
\end{lem}
Let $\psi(x)= e^{xH}$ be the fundamental solution of the linear
system
\[
z'(x)=Hz(x) \qquad x \in [0,1]
\]
where $H$ is a semi-simple infinitesimally symplectic matrix and let
us denote by $L_0'$ the vertical Lagrangian $\oplus_{j=1}^p L_0^j$
of the symplectic space $(V, \oplus_{j=1}^{p}\omega_j)$ where
$\omega_j$ is the standard symplectic form in $\R^{2m}$ for $m =1,2$
corresponding to the decomposition of $V$ into $2$ and $4$
dimensional $\psi(1)$-invariant symplectic subspaces. Then as a
direct consequence of the product property of the Maslov index the
following holds.
\begin{prop}\label{thm:Teoremafinalemaslovnormaleloop}
Let $e^{i\alpha_1}, \dots, e^{i\alpha_k}$ be the Krein positive
eigenvalues of $\psi$. Then the Maslov index with respect to the
Lagrangian $L_0'$ is given by:
\begin{equation}\label{eq:formuladelmaslovnormaleperloop}
\mu_{L'_0}(\psi)\,=\, \sum_{j=1}^{k}f\Big(\frac{\alpha_j}{\pi}\Big),
\end{equation}
where $f$ be the function which holds identity on semi-integer and
is the closest semi-integer not integer otherwise.
\end{prop}
\begin{remark}
By using the zero property for the Maslov index (see for instance
\cite{RobSal}) a direct computation shows that the (ii) and (iii) of
Definition \ref{def: normalisemisemplicisimplectiic3} do not give
any non null contribution to the Maslov index.
\end{remark}

\subsection{The Kashiwara and H\"ormander index} The aim of this
section is to discuss a different notion of Maslov index. Our basic
references are  \cite{Dui76}, \cite{LioVer},
\cite[Section~8]{CapLeeMil} and \cite[Section 3]{GiaPicPor}. The
H\"ormander index, or four-fold index has been introduced in
\cite[Chapter 10, Sect. 3.3]{Hor} who also gave an explicit formula
in terms of a triple of Lagrangian subspaces, which is known in
literature with the name of {\em Kashiwara index\/} and which we now
describe.

Given three Lagrangians $L_1,L_2,L_3\in\Lambda(V,\omega)$ the {\em
Kashiwara index\/} $\tau_V(L_1,L_2,L_3)$ is defined as the signature
of the (symmetric bilinear form associated to the) quadratic form
$Q:L_1\oplus L_2\oplus L_3\to\R$ given by:
\begin{equation}\label{eq:formulapercalcolodiKashiwara}
Q(x_1,x_2,x_3)=\omega(x_1,x_2)+\omega(x_2,x_3)+\omega(x_3,x_1).
\end{equation}
It is proven in \cite[Section~8]{CapLeeMil} that $\tau_V$ is the
unique integer valued map on $\Lambda\times\Lambda\times\Lambda$
satisfying the following properties:
\begin{itemize}
\item[\textbf{[P1]}] (skew symmetry) If $\sigma$ is a permutation of
the set $\{1,2,3\}$,
\[\tau_V(L_{\sigma(1)},L_{\sigma(2)},L_{\sigma(3)})=\textrm{sign}(\sigma)
\,\tau_V(L_1,L_2,L_3);\]
\item[\textbf{[P2]}] (symplectic additivity) given symplectic spaces
$(V,\omega)$, $(\widetilde V,\widetilde\omega)$, and Lagrangians
$L_1,L_2,L_3\in\Lambda(V,\omega)$, $\widetilde L_1, \widetilde L_2,
\widetilde L_3 \in\Lambda(\widetilde V,\widetilde \omega)$, then:
\[\tau_{V\oplus \widetilde V}(L_1\oplus \widetilde L_1,L_2\oplus \widetilde L_2,
L_3\oplus \widetilde L_3)=\tau_V(L_1,L_2,L_3)+\tau_{\widetilde
V}(\widetilde L_1, \widetilde L_2,\widetilde L_3);\]
\item[\textbf{[P3]}] (symplectic invariance) if $\phi:(V,\omega)\to(\widetilde
V,\widetilde\omega)$ is a symplectomorphism, then:
\[\tau_V(L_1,L_2,L_3)=\tau_{\widetilde V}(\phi(L_1),\phi(L_2),\phi(L_3));\]
\item[\textbf{[P4]}] (normalization) if $V=\R^2$ is endowed with
the canonical symplectic form, and $L_1=\R(1,0)$, $L_2=\R(1,1)$,
$L_3=\R(0,1)$, then\[\tau_V(L_1,L_2,L_3)=1.\]
\end{itemize}
Let $(V, \omega)$ be a $2n$-dimensional symplectic vector space and
let $L_1, L_2, L_3$ be three Lagrangians and let us assume that
$L_3$ is transversal both to $L_1$ and $L_2$. If $L_1$ and $L_2$ are
transversal we can choose  coordinates $z =(x,y) \in V$ in such a
way that $L_1$ is defined by the equation $y=0$, $L_2$ by the
equation $x=0$ and consequently $L_3$ is defined by $ y=Ax$ for some
symmetric non-singular matrix $A$. We claim that
\[
\tau_V(L_1, L_2, L_3)= \sgn\,A.
\]
In fact, taking into account that every symplectic vector space of
dimension $2n$ is symplectically isomorphic to $(\R^{2n},
\omega_0)$, property [P3] on the symplectic invariant of the
Kashiwara index and equation
\eqref{eq:formulapercalcolodiKashiwara}, it is enough to compute the
signature of the quadratic form $Q$ where
\[ x_1=(x,0), \ \ x_2=(z, Az), \ \ x_3=(0, y),
\ \ \textrm{for}\ \  x,y,z \in \R^n,\] and $\omega =\omega_0$.
Therefore, $\omega_0(x_1, x_2)= \langle z, Az \rangle$,
$\omega_0(x_3, x_1)= -\langle x, y \rangle$ and finally
$\omega_0(x_2, x_3)= \langle x, y \rangle$ and by this we conclude
that
\[
Q(x_1, x_2, x_3)= \langle z, Az \rangle.
\]
In the general case, let $K = L_1 \cap L_2$. Then $K$ is an
isotropic linear subspace of the symplectic space $(V,\omega)$ and
$K^\#/K:=V^K$ is a symplectic vector space with the symplectic form
induced by $(V, \omega)$. If $L$ is any Lagrangian subspace in $(V,
\omega)$ then $L^K = L \cap K^\#  \mod K$, is a Lagrangian subspace
in $V^K$.
\begin{lem}\label{thm:checipossiamoridurre}
For an arbitrary subspace $K$ of $L_1 \cap L_2 + L_2 \cap L_3 + L_3
\cap L_1$,
\[\tau_V(L_1, L_2, L_3)=\tau_{V^K}(L^K_1, L^K_2, L^K_3).\] where
for $i = 1,2,3$ the Lagrangian subspaces $L_i^K$ are the image of
$L_i$ under the symplectic reduction
\[
(K+ K^{\#}) \to V^K :=(K+ K^{\#})/(K\cap K^{\#}).
\]
\end{lem}
\proof For the proof of this result see \cite[Proposition
1.5.10]{LioVer}.\finedim

We will now proceed to a geometrical description of $\tau_V$ using
the Maslov index for paths; to this aim we will introduce the
 {\em H\"ormander index\/}.
\begin{lem}\label{thm:perdefq}
Given four Lagrangians $L_0,L_1,L_0',L_1'\in\Lambda$ and any
continuous curve $l:[a,b]\to\Lambda$ such that $l(a)=L_0'$ and
$l(b)=L_1'$, then the value of the quantity $\mu(l,L_1)-\mu(l, L_0)$
does {\em not\/} depend on the choice of $l$.
\end{lem}
\begin{proof}
See \cite[Theorem 3.5]{RobSal}
\end{proof}
We are now entitled to define the map $
s:\Lambda\times\Lambda\times\Lambda\times\Lambda\longrightarrow
\frac12\Z.$
\begin{defin}\label{thm:deffourfold}
Given $L_0,L_1,L_0',L_1'\in\Lambda$, the {\em H\"ormander index\/}
$s(L_0,L_1;L_0',L_1')$ is the half-integer $\mu(l,L_1)-\mu(l, L_0)$,
where $l:[a,b]\to\Lambda$ is any continuous curve joining
$l(a)=L_0'$ with $l(b)=L_1'$.
\end{defin}
The  H\"ormander's index, satisfies the following symmetries. (See,
for instance \cite[Proposition 3.23]{GiaPicPor}).
We can now establish the relation between the H\"ormander index $s$
and the Kashiwara index $\tau_V$.This will be made in the same way
as in \cite[Section 3]{GiaPicPor}. We define
$\overline{s}:\Lambda\times\Lambda\times\Lambda\to\Z$ by:
\begin{equation}\label{eq:defqbarra}
\overline{s}(L_0,L_1,L_2):=2s(L_0,L_1;L_2,L_0).
\end{equation}
Observe that the function $s$ is completely determined by
$\overline{s}$, because of the following identity:
\begin{multline}\label{eq:relqqbar}2s(L_0,L_1;L_0',L_1')=2s(L_0,L_1;L_0',L_0)+
2s(L_0,L_1;L_0,L_1')\\
=\overline{s}(L_0,L_1,L_0')-\overline{s}
(L_0,L_1,L_1').\end{multline}
\begin{prop}\label{thm:q=tauV}
The map $\overline{s}$ defined in \eqref{eq:defqbarra} coincides
with the Kashiwara index $\tau_V$.
\end{prop}
\begin{proof}
By  uniqueness, it suffices to prove that $\overline{s}$ satisfies
the properties [P1], [P2], [P3] and [P4]. See \cite{GiaPicPor} for
further details.\finedim

\end{proof}
As a direct consequence of Proposition \ref{thm:q=tauV} and formula
\eqref{eq:relqqbar} we have:
\begin{equation}\label{eq:relstau}
s(L_0,L_1;L_0',L_1')=\frac12[\tau_V(L_0,L_1,L_0')-\tau_V
(L_0,L_1,L_1')].
\end{equation}

\section{The main result}\label{sec:Maslovconleyzehnder}
Let $\psi$ be the fundamental solution of the linear Hamiltonian
system
\[
w'(x)= H w(x), \qquad x \in [0,1].
\]
By Lemma \ref{thm:lagrangianiprodotto} we have  $\mu_{L_0}(\psi) =
\mu(\Gr(\psi), L_0 \times L_0)$; hence
\begin{eqnarray}\label{eq: passaggiopsifi}
\mu(\Gr(\psi), L_0 \times L_0)&=&\mu(\Gr(\psi), \Delta) +
s(\Delta,L_0 \times L_0, \Gr(I), \Gr(\psi(1)))=\\\nonumber &=&
\mu_{CZ}(\psi)-\frac12\tau_V(\Delta, L_0 \times  L_0, \Gr(\psi(1))).
\end{eqnarray}
For one periodic loop the last term in formula \eqref{eq:
passaggiopsifi} vanishes identically because of the anti-symmetry of
the Kashiwara index and by the fact that $\Gr(\psi(1))= \Delta$.
Thus in this case we conclude that
\[
\mu_{L_0}(\psi)\, = \,\mu_{CZ}(\psi).
\]
From now on we assume the following transversality condition:
\begin{enumerate}
\item[(H)] $\psi(1)L_0 \cap L_0 = \{0\}.$
\end{enumerate}
Let $L= L_0 \times L_0$ and $L_2= \Gr\,(\psi(1))$. Thus we only need
to compute the last term in formula \eqref{eq: passaggiopsifi} which
is $-\frac12\tau_V(\Delta, L, L_2)$ where the product form can be
represented with respect to the scalar product in $\R^{4n}$ by the
matrix
\begin{equation*}
\widetilde J=\begin{pmatrix} -J & 0 \\
0 & J
\end{pmatrix}.
\end{equation*} for $J$ defined in \eqref{eq:strutcomplessa}.
We denote by $K$ the isotropic subspace $\Delta \cap L$; it is the
set of all vectors of the form $(0, u, 0, u)$ for $u \in \R^n$.
Moreover $K^{\#}$ is
\[ K^{\#}=\{(x, y, z, v)\in \R^{4n}| \ \ \overline \omega[(x, y, z, v), (0,u,0,u)^T] =0 \}
=\{(x, y, x, v) \colon x, y, v \in \R^n\}.
\]
Identifying the quotient space $K^\#/K$ with the orthogonal
complement $S_K$ of $K$ in $K^\#$ we have $ S_K= \{(t, w, t, -w)
\colon t, w \in \R^n\};$ moreover $\Delta \cap K^{\#}=\Delta$,
$L\cap K^{\#}= L$. Now if $\psi(1)$ has the following block
decomposition
\[\psi(1)=
\begin{pmatrix}
A&B\\C&D
\end{pmatrix},
\]
then  $ L_2 \cap K^{\#}= \{\big(r, s ,Ar + Bs, Cr + Ds\big)\colon
Ar+Bs=r;  r , s \in \R^n\}.$ Since $K^{\#} = S_K \oplus K$ then the
image in $S_K$ of an arbitrary point in $K^\#$ is represented by the
point $[\epsilon, \eta, \epsilon, - \eta]$ where $[\cdot]$ denotes
the equivalence class in the quotient space. Thus we have
\begin{itemize}
\item[(i)] $\Delta^K = \{[\alpha,0,\alpha,0]; \alpha \in \R^n\}$;
\item[(ii)]  $L^K=\{[0,u,0,-u]; u \in \R^n\}$;
\item[(iii)] $L_2^K=\{[r,s,r, Cr+Ds]; Ar + Bs =r\} =  \{[r,s-Ds,r, Cr+Ds]; Ar + Bs =r\}$.
\end{itemize}
Then we have:
\begin{eqnarray*}
\bar \omega(x_1, x_2)&=&\bar
\omega\big([\alpha,0,\alpha,0],[0,u,0,-u]
\big)= -2\langle \alpha, u\rangle.\\
\bar \omega(x_2,x_3)&=&\bar\omega\left([0,u,0,-u],
\left[r,s-Ds,r,Cr\right]\right)= 2\langle u,
r\rangle.\\
\bar \omega(x_3, x_1)&=& \bar \omega\left(\left[r,s-Ds,r,Cr\right],
[\alpha,0,\alpha,0]\right)=\langle s-Ds-Cr, \alpha \rangle .
\end{eqnarray*}
Hence the quadratic form $Q$ is given by:
\[
Q(x_1,x_2, x_3)= -2\langle\alpha, u\rangle+2\langle u,
r\rangle+\langle s-Ds-Cr,  \alpha\rangle,\ \ \   \alpha, u, r, s,
\in \R^n \ \textrm{and} \ Ar + Bs =r.
\]
Due to the transversality condition $(H)$ we have $s= B^{-1}(I_n-A)
r$ and by setting $2X = (I_n-D)B^{-1}(I_n -A) -C$ the quadratic form
$Q$ can be written as follows
\[
Q(x_1,x_2, x_3)= -2\langle\alpha, u\rangle+2\langle u, r\rangle+ 2
\langle X r, \alpha\rangle = \langle Y w, w \rangle
\]
for $w =(\alpha, u, r)$ and $Y$ given by the matric below:
\[
Y= \begin{pmatrix} 0_n & - I_n & X\\
-I_n & 0_n & I_n\\
X^T & I_n & 0_n
\end{pmatrix}.
\]
The Cayley-Hamilton polynomial of $A$ is given by
\[
p_Y(\lambda)= \lambda^3 I_n -(2I_n + XX^T)\lambda + (X + X^T).
\]
\begin{lem}\label{lem:lemmafava} Given the symplectic block matrix $
\begin{pmatrix}
A&B\\C&D
\end{pmatrix}$
then the $n$ by $n$ matrix $2X = (I_n-D)B^{-1}(I_n -A) -C$ is
symmetric.
\end{lem}
\proof In fact
\begin{eqnarray*}
2X &=& (B^{-1}-DB^{-1})(I-A)-C= B^{-1} - B^{-1}A - DB^{-1} + DB^{-1}A -C \\
2X^T &=&
[B^T]^{-1}-[B^T]^{-1}D^T - A^T[B^T]^{-1}+ A^T[B^T]^{-1}D^T-C^T.
\end{eqnarray*}
Moreover by multiplying this last equation on the right by $BB^{-1}$
we have
\begin{eqnarray*}
2X^T
&=&([B^T]^{-1}B-[B^T]^{-1}D^TB-A^T[B^T]^{-1}B+A^T[B^T]^{-1}D^TB-C^TB)B^{-1}=
\\&=&([B^T]^{-1}B-[B^T]^{-1}B^TD-A^T[B^T]^{-1}B+A^T[B^T]^{-1}B^TD-C^TB)B^{-1}=\\
&=&([B^T]^{-1}B-D-A^T[B^T]^{-1}B+A^TD-C^TB)B^{-1}=\\
&=&([B^T]^{-1}B-D-A^T[B^T]^{-1}B+I_n)B^{-1}=[B^T]^{-1}-DB^{-1}-
A^T[B^T]^{-1}+ B^{-1},
\end{eqnarray*}
where we used the relations $A^TD-C^TB=I_n$, $A^TC=C^TA$ and finally
$D^TB=B^TD$. Thus by the expression for $2X$ and this last equality
we conclude that in order to prove the thesis it is enough to show
that:
\[
-B^{-1}A + DB^{-1}A -C =[B^T]^{-1}-A^T[B^T]^{-1}.
\]
Now we observe that by multiplying on the left the first member of
the above equality by $[B^T]^{-1}B^T$, we have
\begin{eqnarray*}
([B^T]^{-1}B^T)(-B^{-1}A + DB^{-1}A -C)&=& [B^T]^{-1}(-B^TB^{-1}A + B^TDB^{-1}A -B^TC)=\\
&=&[B^T]^{-1}(-B^TB^{-1}A + D^TBB^{-1}A -B^TC)=\\
&=&[B^T]^{-1}(-B^TB^{-1}A + D^TA -B^TC)=\\
&=&-B^{-1}A+[B^T]^{-1}.
\end{eqnarray*}
Thus we reduced to show that
$-B^{-1}A+[B^T]^{-1}=[B^T]^{-1}-A^T[B^T]^{-1}$ or which is the same
to $B^{-1}A=A^T[B^T]^{-1}$. Otherwise stated since $A^T[B^T]^{-1}=
(B^{-1}A)^T$, it is enough to check that the $n$ by $n$ matrix
$U=B^{-1}A$ is symmetric. In fact
\[
U=I_n\cdot U = (A^TD-C^TB)B^{-1}A=A^TDB^{-1}A-C^TA
=A^TDB^{-1}A-A^TC= A^T(DB^{-1}A-C);
\]
moreover $U^T=A^T[B^T]^{-1}$. Thus the condition $U^T=U$ reduced to
$A^T[B^T]^{-1}=A^T(DB^{-1}A-C)$ and then the only thing to show is
that $[B^T]^{-1}=DB^{-1}A-C$. In fact by multiplying the second
member on the left by $[B^T]^{-1}B^T$, it then follows:
\[
[B^T]^{-1}(B^TDB^{-1}A-B^TC)=[B^T]^{-1}(D^TBB^{-1}A-B^TC)=[B^T]^{-1}(D^TA-B^TC)=[B^T]^{-1}I_n
\]
and this conclude the proof .\finedim

Since $\psi(1)$ is a symplectomorphism and $X$ is a symmetric matrix
by Lemma \ref{lem:lemmafava}, there exists an $n$ by $n$ orthogonal
matrix such that $M^T X M = \diag(\lambda_1, \dots, \lambda_k)$
where the eigenvalues $\lambda_j$ are counted with multiplicity;
hence in order to compute the solutions of $p_Y(\lambda)=0$ it is
enough to compute the solutions of $M^T p_Y(\lambda) M=0$ which is
the same to solve
\[
0=\lambda^3  -(2 + \lambda_j^2)\lambda +
2\lambda_j=(\lambda-\lambda_j)(\lambda^2 + \lambda_j \lambda -2), \
\ \textrm{for}\ \  j=1, \dots, k.
\]
The solutions of the equation $\lambda^2 + \lambda_j \lambda -2=0$
are one positive and one negative and therefore they do not give any
contribution in the signature of $Y$. Thus we proved that $
\sgn\,(Y) = \sgn\,(X)$ and therefore
\[
-\frac12\tau_{\R^{4n}}(\Delta, L, L_2)=-\frac12
\,\tau_{\R^{2n}}(\Delta^K, L^K, L_2^K)= -\frac12\, \sgn\,X.
\]
Summing up the previous calculation we proved the following result.
\begin{thm}\label{thm:Teoremafinalemaslov}
Let $\psi\colon [0,1] \to \Spl(2n)$ be the fundamental solution of
the Hamiltonian system
\[
w'(x)= Hw(x), \qquad x \in [0,1]
\]
and let us assume that condition $(H)$ holds.
 Then the Maslov index of $\psi$ is given by
\begin{equation}\label{eq:finale}
\mu_{L_0}(\psi)= \mu_{CZ}(\psi) +\frac12\,\sgn\, \widetilde X
\end{equation}
for $\widetilde X = C+ (D-I_n)B^{-1}(I_n -A) $. In particular if
$\psi(1)=I_{2n}$ formula \eqref{eq:finale} reduces to
\[
\mu_{L_0}(\psi)= \mu_{CZ}(\psi).
\]
\end{thm}
\begin{cor} \label{coro2.3} Let $e^{i \alpha_1}, \dots, e^{i
\alpha_k}$ be the Krein positive purely imaginary eigenvalues of the
fundamental solution $\psi(x)=e^{xH}$ counted with algebraic
multiplicity and we assume that $(H)$ and $\det(\psi(1) -
I_{2n})\not=0$ hold. Then the Maslov index of $\psi$ is given by
\begin{equation}
\mu_{L_0}(\psi)= \sum_{j=1}^{k}g\Big(\frac{\alpha_j}{\pi}\Big)+
\frac12\,\mathop{\rm sign} \widetilde X,
\end{equation}
where we denoted by $g$ the double integer part function which holds
the identity on integers and it is the closest odd integer
otherwise.
\end{cor}

\end{document}